\RequirePackage{etoolbox}
\csdef{input@path}{{style/}{graphics/}}
\documentclass[aap,MSNbibl,dvips]{arximspdf}
\usepackage{url,breakurl}

\doi{10.1214/14-AAP1002} %kopijuoti is PTS
\volume{25}
\issue{2}
\pubyear{2015}
\firstpage{523}
\lastpage{547}

\makeatletter
\newcommand{\eqref}[1]{(\ref{#1})}
\newtheorem{thmm}{Theorem}
\newtheorem{lem}[thmm]{Lemma}

\newproclaim{definition}[thmm]{Definition}
\newproclaim{ex}[thmm]{Example}
\newproclaim{rem}[thmm]{Remark}

\newproclaim{remu}{Remark}

\newproclaim{continua}{Example~\ref{esem1} (continued)}
\newproclaim{cont}{Example~\ref{eped} (continued)}

\makeatother

\begin{document}
\begin{frontmatter}
\title{Central limit theorems for an Indian buffet model with random weights}
\runtitle{An Indian buffet model}

\begin{aug}
\author[A]{\fnms{Patrizia} \snm{Berti}\ead[label=e1]{patrizia.berti@unimore.it}},
\author[B]{\fnms{Irene} \snm{Crimaldi}\thanksref{T1}\ead[label=e2]{irene.crimaldi@imtlucca.it}},
\author[C]{\fnms{Luca} \snm{Pratelli}\ead[label=e3]{pratel@mail.dm.unipi.it}}
\and
\author[D]{\fnms{Pietro} \snm{Rigo}\corref{}\ead[label=e4]{pietro.rigo@unipv.it}}
\runauthor{Berti, Crimaldi, Pratelli and Rigo}
\thankstext{T1}{Supported in part by
CNR PNR project ``CRISIS Lab.''}
\affiliation{Universit\`{a} di Modena e Reggio-Emilia, IMT Institute
for Advanced Studies, Accademia Navale and Universit\`{a} di Pavia}
\address[A]{P. Berti\\
Dipartimento di Matematica Pura\\
\quad ed Applicata ``G. Vitali''\\
Universit\`{a} di Modena e Reggio-Emilia\\
via Campi 213/B\\
41100 Modena\\
Italy\\
\printead{e1}} %adresu isvedimo komanda gale!
\address[B]{I. Crimaldi\\
IMT Institute for Advanced Studies\hspace*{22pt}\\
Piazza San Ponziano 6\\
55100 Lucca\\
Italy\\
\printead{e2}}
\address[C]{L. Pratelli\\
Accademia Navale\\
viale Italia 72\\
57100 Livorno\\
Italy\\
\printead{e3}}
\address[D]{P. Rigo\\
Dipartimento di Matematica ``F. Casorati''\\
Universit\`{a} di Pavia\\
via Ferrata 1\\
7100 Pavia\\
Italy\\
\printead{e4}}
\end{aug}

% HISTORY:
\received{\smonth{4} \syear{2013}}
\revised{\smonth{10} \syear{2013}}

% ABSTRACT
%
\begin{abstract}
The three-parameter Indian buffet process is generalized. The possibly
different role played by customers is taken into account by suitable
(random) weights. Various limit theorems are also proved for such
generalized Indian buffet process. Let $L_n$ be the number of dishes
experimented by the first $n$ customers, and let $\overline{K}_n=(1/n)
\sum_{i=1}^nK_i$ where $K_i$ is the number of dishes tried by customer
$i$. The asymptotic distributions of $L_n$ and $\overline{K}_n$,
suitably centered and scaled, are obtained. The convergence turns out
to be stable (and not only in distribution). As a particular case, the
results apply to the standard (i.e., nongeneralized) Indian buffet process.
\end{abstract}

% KEYWORDS

\end{frontmatter}
%
%s1 #&#
\section{Introduction}\label{intro}

Let $(\mathcal{X},\mathcal{B})$ be a measurable space. Think of
$\mathcal{X}$ as a collection of features potentially shared by an
object. Such an object is assumed to have a finite number of features
only and is identified with the features it possesses. To investigate
the object, thus, we focus on the finite subsets of $\mathcal{X}$.

Each finite subset $B\subset\mathcal{X}$ can be associated to the
measure $\mu_B=\sum_{x\in B}\delta_x$, where $\mu_\varnothing=0$ and
$\delta_x$ denotes the point mass at $x$. If $B$ is random, $\mu_B$
is random as well. In fact, letting $F=\{\mu_B\dvtx B$ finite$\}$, there
is a growing literature focusing on those random measures $M$
satisfying $M\in F$ a.s. See \cite{BJP} and most references quoted
below in this section.

A remarkable example is the \emph{Indian Buffet Process} (IBP)
introduced by Griffiths and Ghahramani and developed by Thibaux and
Jordan; see \cite{GG06,GG11,TJ}. The objects are the
customers which sequentially enter an infinite buffet $\mathcal{X}$ and
the features are the dishes tasted by each customer. In this framework,
each customer is modeled by a (completely) random measure $M$ such that
$M\in F$ a.s. The atoms of $M$ represent the dishes experimented by the
customer.

Our starting point is a three-parameter extension of IBP, referred to
as \emph{standard} IBP in the sequel, introduced in \cite{BJP} and
\cite
{TG} to obtain power-law
behavior. Fix \mbox{$\alpha>0$}, $\beta\in[0,1)$ and $c>-\beta$. Here,
$\alpha
$ is the mass parameter, $\beta$ the discount parameter (or stability
exponent) and $c$ the concentration parameter. Also, let $\operatorname{Poi}(\lambda)$
denote the Poisson distribution with mean $\lambda\geq0$, where
$\operatorname{Poi}(0)=\delta_0$. The dynamics of the standard IBP is as follows.
Customer 1 tries $\operatorname{Poi}(\alpha)$ dishes. For each $n\geq1$, let $S_n$ be
the collection of dishes experimented by the first $n$ customers. Then:

\begin{itemize}[$-$]
\item[$-$] Customer $n+1$ selects a subset $S_n^*\subset S_n$. Each
$x\in S_n$ is included or not into $S_n^*$ independently of the other
members of $S_n$. The inclusion probability is
\[
\frac{\sum_{i=1}^nM_i\{x\}-\beta}{c+n},
\]
where $M_i\{x\}$ is the indicator of the event $\{$customer $i$ selects
dish $x\}$.

\item[$-$] In addition to $S_n^*$, customer $n+1$ also tries
$\operatorname{Poi}(\lambda_n)$ new dishes, where $\lambda_n=\alpha\frac{\Gamma
(c+1)\Gamma(c+\beta+n)}{\Gamma(c+\beta)\Gamma(c+1+n)}$.
\end{itemize}

For $\beta=0$, such a model reduces to the original IBP of
\cite{GG06,GG11,TJ}.

IBP is a flexible tool, able to capture the dynamics of various real
problems. In addition, IBP is a basic model in Bayesian
nonparametrics; see \cite{GB} and \cite{KG}. In factor analysis, for
instance, IBP works as an infinite-capacity prior over the space of
latent factors; see \cite{KG}. In this way, the number of factors is
not specified in advance but is inferred from the data. Such a
number is also allowed to grow as new data points are observed.
Among the other possible applications of IBP, we mention causal
inference \cite{WGG}, modeling of choices \cite{GJR}, similarity
judgements~\cite{NG} and dyadic data \cite{MGNR}.

Despite its prominent role, however, the asymptotics of IBP is largely
neglected. To the best of our knowledge, the only known fact is the
a.s. behavior of $L_n$ (defined below) and some other related quantities
for large $n$; see \cite{BJP} and \cite{TG}. Nothing is known as
regards limiting distributions.

In this paper, we aim to do two things:

First, we generalize the standard IBP. Indeed, the discount parameter
$\beta$ is allowed to take values in $(-\infty, 1)$ rather than in
$[0,1)$. More importantly, the possible different relevance of
customers is taken into account by random weights. Let $R_n>0$ be the
weight attached to customer $n$. Then, for each $x\in S_n$, the
inclusion probability becomes
\[
\frac{\sum_{i=1}^nR_iM_i\{x\}-\beta}{c+\sum_{i=1}^nR_i}.
\]
Similarly, the new dishes tried by customer $n+1$ are now
$\operatorname{Poi}(\Lambda_n)$ rather than $\operatorname{Poi}(\lambda_n)$, where
$\Lambda_n=\alpha\frac{\Gamma(c+1)\Gamma(c+\beta+\sum_{i=1}^nR_i)}{\Gamma(c+\beta)\Gamma(c+1+\sum_{i=1}^nR_i)}$.
If $\beta\in[0,1)$ and $R_n=1$ for all $n$, the model reduces to the
standard IBP.

Second, we investigate the asymptotics of the previous generalized IBP
model. We focus on
\begin{eqnarray*}
L_n&=& \mbox{number of dishes experimented by the first }n\mbox{
customers\quad and}
\\
\overline{K}_n&=&\frac{1}{n}\sum_{i=1}^nK_i\qquad
\mbox{where } K_i= \mbox {number of dishes tried by customer }i.
\end{eqnarray*}
Three results are obtained. Define $a_n(\beta)=\log{n}$ if $\beta=0$
and $a_n(\beta)=n^\beta$ if $\beta\in(0,1)$. Then, under some
conditions on the weights $R_n$ (see Theorems \ref{thmln}, \ref{abg},
\ref{cert5r}) it is shown that:

\begin{longlist}[(iii)]
\item[(i)] if $\beta\in[0,1)$, then $\frac{L_n}{a_n(\beta
)}\stackrel
{\mathrm{a.s.}}\longrightarrow\lambda$ where $\lambda>0$ is a certain constant;

\item[(ii)] if $\beta\in[0,1)$, then $\sqrt{a_n(\beta)}  \{
\frac
{L_n}{a_n(\beta)}-\lambda \}\longrightarrow\mathcal{N} (0,
\lambda )$ stably;

\item[(iii)] if $\beta<1/2$, then $\overline{K}_n\stackrel
{\mathrm{a.s.}}\longrightarrow Z$ and
\begin{eqnarray*}
\sqrt{n} \{\overline{K}_n-Z \}&\longrightarrow&\mathcal {N} \bigl(0,
\sigma^2 \bigr)\qquad \mbox{stably},
\\
\sqrt{n} \bigl\{\overline{K}_n-E (K_{n+1}\mid\mathcal
{F}_n ) \bigr\}&\longrightarrow&\mathcal{N} \bigl(0, \tau^2
\bigr) \qquad\mbox{stably},
\end{eqnarray*}
where $Z$, $\sigma^2$, $\tau^2$ are suitable \emph{random} variables,
and $\mathcal{F}_n$ is the sub-$\sigma$-field induced by the
available information at time $n$.
\end{longlist}

Stable convergence is a strong form of convergence in
distribution. The basic definition is recalled in Section~\ref{sc6h}.
Further, $\mathcal{N}(0,a)$ denotes the Gaussian law
with mean 0 and variance $a\geq0$, where
$\mathcal{N}(0,0)=\delta_0$.

Among other things, the above results can be useful in making
(asymptotic) inference on the model. As an example, suppose $\beta\in
[0,1)$. In view of~(i),
\[
\widehat{\beta}_n=\frac{\log{L_n}}{\log{n}}
\]
is a strongly consistent estimator of $\beta$ for each $\beta\in
[0,1)$. In turn, (ii) provides the
limiting distribution of $\widehat{\beta}_n$ so that simple tests on
$\beta$ can be manufactured. Similarly, if $\beta<1/2$, asymptotic
confidence bounds for the random limit $Z$ of $\overline{K}_n$ can be
obtained by (iii); see Section~\ref{v4gr}.

Note also that, because of (iii), the convergence rate of
$\overline{K}_n-E (K_{n+1}\mid\mathcal{F}_n )$ is at least
$n^{-1/2}$. Therefore, $\overline{K}_n$ is a good predictor of
$K_{n+1}$ for large $n$ and $\beta<1/2$; see Section~\ref{v4gr} again.

The results in (i)--(iii) hold in particular if $R_n=1$ for all
$n$. Thus, (ii) and (iii) provide the limiting distributions of
$L_n$ and $\overline{K}_n$ in the standard IBP model. Furthermore, in
this case, (iii) holds for all $\beta<1$ and not only for $\beta<1/2$.\vadjust{\goodbreak}

We close this section with some remarks on $\beta$ and the $R_n$.

\textit{The discount parameter $\beta$.} Roughly speaking, if $\beta<0$,
the inclusion probabilities are larger and the chances of
tasting new dishes vanish very quickly; see Lemma~\ref{basic123l}.
Define in fact
\[
L=\sup_nL_n=\operatorname{card}\{x\in
\mathcal{X}\dvtx x\mbox{ is tried by some customer}\}.
\]
Because of (i), $L_n$ increases logarithmically if
$\beta=0$ while exhibits a power-law behavior if $\beta\in(0,1)$.
Accordingly, $L=\infty$ a.s. if $\beta\in[0,1)$. On the contrary,
\[
E \bigl(e^L \bigr)<\infty\qquad\mbox{if }\beta<0;
\]
see Lemma~\ref{pino78}. In particular, $\beta<0$ implies $L<\infty$
a.s., and this fact can help to describe some real situations.

Formally, the model studied
in this paper makes sense whenever $R_n>\max(\beta,0)$ for all $n$.
Hence, one could also admit $\beta\geq
1$. However, $\beta=1$ leads to
trivialities. Instead, $\beta>1$ could be potentially interesting, but
it is hard to unify the latter case and $\beta<1$. Accordingly, we will
focus on $\beta<1$.

Unless $R_n=1$ for all $n$, the results in (iii) are available for
$\beta<1/2$ only. Certainly, (iii) can fail
if $\beta\in[1/2, 1)$. Perhaps, some form of (iii) holds even if
$\beta\in[1/2, 1)$, up to replacing $\sqrt{n}$ with some other
norming constant and $\mathcal{N} (0, \sigma^2 )$ and
$\mathcal
{N} (0, \tau^2 )$ with some other limit kernels. But we did not
investigate this issue.

A last note is that
$\beta$ plays an analogous role to that of the discount parameter in
the two-parameter Poisson--Dirichlet process. Indeed, such parameter
regulates the asymptotic
behavior of the number of distinct observed values, in the same way
as $\beta$ does for $L_n$. See, for example, \cite{BCL,PY,PIT} for the two-parameter Poisson--Dirichlet and \cite{BJP,TG} for the standard IBP.

\textit{The weights $R_n$.} Standard IBP has been generalized in various
ways, mainly focusing on computational issues; see, for example, \cite
{DG,GFB,MGJ,WOG}. In this paper, the possible
need of distinguishing objects according to some associated random
factor is dealt with. To this end, customer $n$ is attached a random
weight $R_n$. Indeed, it may be that different customers have different
importance, due to some random cause, that does not affect their
choices but is relevant to the choices of future customers. Analogous
models occur in different settings, for instance in connection with P\'
olya urns and species sampling sequences; see \cite{AMS,BCL,BPR04,BCPR09,BCPR11,PEM}.

The model investigated in this paper, referred to as ``weighted'' IBP
in the sequel, generally applies to evolutionary phenomena. In a
biological framework, for instance, a newborn exhibits some features in
common with the existing units with a probability depending on the
latter's weights (reproductive power, ability of adapting to new
environmental conditions or to compete for finite resources, and so
on). The newborn also presents some new features that, in turn, will be
transmitted to future generations with a probability depending on
his/her weight. See, for example, \cite{BFF} and \cite{RICE}.

Similar examples arise in connection with the evolution of language;
see, for example, \cite{CROFT}. A neologism (i.e., a newly coined term,
word, phrase or concept) is often directly attributable to a specific
people (or journal, period, event and so on) and its diffusion depends
on the importance of such a people. For instance, suppose we are given
a sample of journals of the same type (customers) during several years.
Each journal uses words (dishes), some of which have been previously
used while some others are new. A word appearing for the first time in
a journal has a probability of being reused which depends on the
importance of the journal at issue.

Other applications of the weighted IBP could be found in Bayesian
nonparametrics. Standard IBP is widely used as a prior on binary
matrices with a fixed finite number of rows and infinitely many columns
(rows correspond to objects and columns to features). The weighted IBP
can be useful in all those settings where customers arrive \textit
{sequentially}. As an example, some dynamic networks present a
competitive aspect, and not all nodes are equally successful in
acquiring links. Suppose the network evolves in time, a node (customer)
is added at every time step and some links are created with some of the
existing nodes. The different ability of competing for links is modeled
by a weight attached to each node; see for example, \cite{BB}.
Following \cite{MGJ2009} and \cite{SCJ}, each node could be described
by a set of binary features (dishes) and the probability of a link is a
function of the features of the involved nodes. A nonparametric latent
feature model could be assessed at every time step, with the weighted
IBP as a prior on the feature matrix.

A last remark concerns the probability distribution of the sequence
$(M_n)$, where $M_n$ is the random measure corresponding to customer
$n$. Because of the weights, unlike the standard IBP, $(M_n)$ can fail
to be exchangeable. Thus, the usual machinery of Bayesian
nonparametrics cannot be automatically implemented, due to the lack of
exchangeability, and this can create some technical drawbacks. On the
other hand, the exchangeability assumption is often untenable in
applications. In such cases, the weighted IBP is a realistic
alternative to the standard IBP. We finally note that, when $\beta=0$,
$(M_n)$ satisfies a weak form of exchangeability known as conditional
identity in distribution; see Section~\ref{kh8u5rf} and Lemma~\ref{c4p}.

%s2 #&#
\section{Preliminaries}

%s2.1 #&#
\subsection{Basic notation}

Throughout, $\mathcal{X}$ is a separable metric space and $\mathcal{B}$
the Borel $\sigma$-field on $\mathcal{X}$. We let
\[
\mathcal{M}=\{\mu\dvtx\mu\mbox{ is a finite positive measure on }\mathcal {B}\},
\]
and we say that $\mu\in\mathcal{M}$ is \textit{diffuse} in case
$\mu\{x\}
=0$ for all $x\in\mathcal{X}$.

All random variables appearing in this paper, unless otherwise stated,
are defined on a fixed probability space $(\Omega,\mathcal{A},P)$. If
$\mathcal{G}\subset\mathcal{A}$ is a sub-$\sigma$-field, and $X$ and
$Y$ are random variables with values in the same measurable space, we write
\[
X\mid\mathcal{G}\sim Y\mid\mathcal{G}
\]
to mean that $P(X\in A\mid\mathcal{G})=P(Y\in A\mid\mathcal{G})$
a.s. for each measurable set $A$.

%s2.2 #&#
\subsection{Random measures}\label{sh6f}

A \emph{random measure} (r.m.) is a map $M\dvtx\Omega\rightarrow
\mathcal{M}$
such that $\omega\mapsto M(\omega)(B)$ is $\mathcal{A}$-measurable for
each $B\in\mathcal{B}$. In the sequel, we write $M(B)$ to denote the
real random variable $\omega\mapsto M(\omega)(B)$. Similarly, if
$f\dvtx\mathcal{X}\rightarrow\mathbb{R}$ is a bounded measurable function,
$M(f)$ stands for
\[
M(\omega) (f)=\int f(x) M(\omega) (dx).
\]

A \emph{completely} r.m. is an r.m. $M$ such that
$M(B_1),\ldots,M(B_k)$ are independent random variables whenever
$B_1,\ldots,B_k\in\mathcal{B}$ are pairwise disjoint; see \cite{KING}.

Let $\nu\in\mathcal{M}$. A \emph{Poisson r.m. with intensity} $\nu
$ is a
completely r.m. $M$ such that $M(B)\sim \operatorname{Poi} (\nu(B) )$ for
all $B\in\mathcal{B}$. Note that $M(B)=0$ a.s. in case $\nu(B)=0$.
Note also that the intensity $\nu$ has been requested to be a finite
measure (and not a $\sigma$-finite measure as it usually happens).

We refer to \cite{CIN} and \cite{KING} for Poisson r.m.'s. We just note
that a Poisson r.m. with intensity $\nu$ is easily obtained. Since
$\nu
$ has been assumed to be a finite measure, it suffices to let $M=0$ if
$\nu(\mathcal{X})=0$, and otherwise
\[
M=I_{\{N>0\}} \sum_{j=1}^N
\delta_{X_j},
\]
where $(X_j)$ is an i.i.d. sequence of $\mathcal{X}$-valued random
variables with $X_1\sim\nu/\nu(\mathcal{X})$, $N$ is independent of
$(X_j)$ and $N\sim \operatorname{Poi} (\nu(\mathcal{X}) )$.

As in Section~\ref{intro}, let $F=\{\mu_B\dvtx B$ finite$\}$ where
$\mu_\varnothing=0$ and $\mu_B=\sum_{x\in B}\delta_x$. Since
$\mathcal{X}$ is separable metric and $\mathcal{B}$ the Borel
$\sigma$-field, the set $\{M\in F\}$ belongs to $\mathcal{A}$ for
every r.m. $M$. In this paper, we focus on those r.m.'s $M$
satisfying $M\in F$ a.s. If $M$ is a Poisson r.m. with intensity
$\nu$, then $M\in F$ a.s. if and only if $\nu$ is diffuse.
Therefore, another class of r.m.'s is to be introduced.

Each $\nu\in\mathcal{M}$ can be uniquely written as $\nu=\nu_c+\nu_d$,
where $\nu_c$ is diffuse and
\[
\nu_d=\sum_j\gamma_j
\delta_{x_j}
\]
for some $\gamma_j\geq0$ and $x_j\in\mathcal{X}$. (The case $\nu_d=0$
corresponds to $\gamma_j=0$ for all $j$.) Say that $M$ is a \emph
{Bernoulli r.m. with hazard measure} $\nu$, where $\nu\in\mathcal
{M}$, if:

\begin{itemize}
\item$M=M_1+M_2$ with $M_1$ and $M_2$ independent r.m.'s;

\item$M_1$ is a Poisson r.m. with intensity $\nu_c$;

\item$M_2=\sum_jV_j \delta_{x_j}$ where the $V_j$ are independent
indicators satisfying $P(V_j=1)=\gamma_j$.
\end{itemize}

Some (obvious) consequences of the definition are the following:
\begin{itemize}[$-$]
\item[$-$] For each $B\in\mathcal{B}$, $E \{M(B)\}=\nu(B)$ and
\[
E \bigl\{M(B)^2 \bigr\}=\nu(B)+\nu(B)^2-\sum
_{x\in B}\nu\{x\}^2;
\]
\item[$-$] $M=M_1$ a.s. if $\nu=\nu_c$ and $M=M_2$ a.s. if $\nu
=\nu_d$;

\item[$-$] $M$ is a completely r.m.;

\item[$-$] $M\in F$ a.s.
\end{itemize}

We will write
\[
M\sim \operatorname{Be}P(\nu)
\]
to mean that $M$ is a Bernoulli r.m. with hazard measure $\nu$.

%s2.3 #&#
\subsection{Stable convergence}\label{sc6h} Stable convergence is a
strong form of convergence in distribution. We just recall the basic
definition and we refer to \cite{CLP,HH} and references therein
for more information.

An r.m. $K$ such that $K(\omega)(\mathcal{X})=1$, for all $\omega\in
\Omega$, is said to be a \emph{kernel} or a \emph{random
probability measure}. Let $K$ be a kernel and $(X_n)$ a sequence of
$\mathcal{X}$-valued random variables. Say that $X_n$ \emph{converges
stably to} $K$ if
\[
E \bigl\{K(f)\mid H \bigr\}=\lim_nE \bigl
\{f(X_n)\mid H \bigr\}
\]
for all $H\in\mathcal{A}$ with $P(H)>0$ and all bounded continuous
$f\dvtx\mathcal{X}\rightarrow\mathbb{R}$. (Recall that $\mathcal{A}$
denotes the basic $\sigma$-field on $\Omega$.) For $H=\Omega$, stable
convergence trivially
implies convergence in distribution.

%s2.4 #&#
\subsection{Conditionally identically distributed sequences}\label{kh8u5rf}

Let $(X_n\dvtx n\geq1)$ be a sequence of random variables (with values
in any measurable space) adapted to a filtration
$(\mathcal{U}_n\dvtx n\geq0)$. Say that $(X_n)$ is \textit{conditionally
identically distributed} (c.i.d.) with respect to $(\mathcal{U}_n)$
in case
\[
X_k\mid\mathcal{U}_n\sim X_{n+1}\mid
\mathcal{U}_n \qquad\mbox{for all }k>n\geq0.
\]
Roughly speaking this means that, at each time $n\geq0$, the future
observations $(X_k\dvtx k>n)$ are identically distributed given the past
$\mathcal{U}_n$. If $\mathcal{U}_0=\{\varnothing,\Omega\}$ and
$\mathcal{U}_n=\sigma(X_1,\ldots,X_n)$, the filtration
$(\mathcal{U}_n)$ is not mentioned at all and $(X_n)$ is just called
c.i.d. Note that $X_k\sim X_1$ for all $k\geq1$ whenever $(X_n)$ is
c.i.d.

The c.i.d. property is connected to exchangeability. Indeed, $(X_n)$
is exchangeable if and only if it is stationary and c.i.d., and the
asymptotic behavior of c.i.d. sequences is quite close to that of
exchangeable ones. We refer to \cite{BPR04} for details.

%s3 #&#
\section{The model}\label{nm}

Let $(M_n\dvtx n\geq1)$ be a sequence of r.m.'s and $(R_n\dvtx n\geq
1)$ a~sequence of real random variables. The probability distribution of
$((M_n,R_n)\dvtx n\geq1)$ is identified by the parameters $m, \alpha,
\beta$ and $c$ as follows:

\begin{itemize}
\item$m$ is a diffuse probability measure on $\mathcal{B}$;

\item$\alpha, \beta, c$ are real numbers such that $\alpha>0$,
$\beta
<1$ and $c>-\beta$;

\item$R_n$ independent of $(M_1,\ldots,M_n,R_1,\ldots,R_{n-1})$ and
$R_n\geq u>\max(\beta, 0)$, for some constant $u$ and each $n\geq1$;

\item$M_{n+1}\mid\mathcal{F}_n\sim \operatorname{Be}P(\nu_n)$ for all $n\geq0$, where
\begin{eqnarray*}
\mathcal{F}_0&=&\{\varnothing,\Omega\},\qquad \nu_0=\alpha m,\qquad
\mathcal {F}_n=\sigma(M_1,\ldots,M_n,R_1,
\ldots,R_n),
\\
\nu_n&=&\sum_{x\in S_n} \frac{\sum_{i=1}^nR_iM_i\{x\}-\beta}{\sum_{i=1}^nR_i+c}
\delta_x + \frac{\Gamma(c+1)\Gamma(c+\beta+\sum_{i=1}^nR_i)}{\Gamma(c+\beta)\Gamma(c+1+\sum_{i=1}^nR_i)} \alpha m
\quad\mbox{and}\\
 S_n&=& \bigl\{x\in\mathcal{X}\dvtx M_i\{x\}=1\mbox{ for some }i=1,\ldots,n \bigr\}.
\end{eqnarray*}
\end{itemize}

Our model is the sequence $((M_n,R_n)\dvtx n\geq1)$. It reduces to the
standard IBP in case $\beta\in[0,1)$ and $R_n=1$ for all $n$. Note
that $M_1$ is a Poisson r.m. with intensity~$\alpha m$. Note also
that $M_n\in F$ a.s. for all $n\geq1$, so that
\[
S_n=\bigcup_{i=1}^n
\mbox{Support}(M_i)\qquad \mbox{a.s.}
\]

Formally, for such a model to make sense, $\beta$ can be taken to be
any real number satisfying $R_n>\max(\beta, 0)$ for all $n$. For the reasons
explained in Section~\ref{intro}, however, in this paper we focus on
$\beta<1$. We also assume $R_n\geq u$, for all $n$ and some constant
$u>\max(\beta, 0)$, as a mere technical assumption. In the sequel, we let
\[
\Lambda_0=\alpha\quad\mbox{and}\quad \Lambda_n=\alpha
\frac{\Gamma
(c+1)\Gamma(c+\beta+\sum_{i=1}^nR_i)}{\Gamma(c+\beta)\Gamma
(c+1+\sum_{i=1}^nR_i)}.
\]
In this notation, the diffuse part of $\nu_n$ can be written as
$\Lambda_n m$.

As remarked in Section~\ref{intro}, $R_n$ should be regarded as the
weight of customer~$n$. Thus, the possibly different role played by
each customer can be taken into account.

Apart from the possible negative values of $\beta$, the parameters
$m, \alpha, \beta$ and $c$ have essentially the same meaning as in
the standard IBP. The probability measure $m$ allows us to draw, at
each step $n\geq1$, an i.i.d. sample of new dishes. In fact,
$m (\mathcal{X}\setminus S_n )=1$ a.s. for $m$ is diffuse
and $S_n$ finite a.s. The mass parameter $\alpha$ controls the total
number of tried dishes per customer. The concentration parameter $c$
tunes the
number of customers which try each dish. The discount parameter
$\beta$ has been discussed in Section~\ref{intro}.

An r.m. can be seen as a random variable with values in
$(\mathcal{M},\Sigma)$, where $\Sigma$ is the $\sigma$-field on
$\mathcal{M}$ generated by the maps $\mu\mapsto\mu(B)$ for all
$B\in\mathcal{B}$. In the standard IBP case, $(M_n)$ is an
exchangeable sequence of random variables. Now, because of the
$R_n$, exchangeability is generally lost. In fact, the same
phenomenon (loss of exchangeability) occurs in various other
extensions of IBP; see \cite{DG,GFB,MGJ,WOG}. However, under some conditions, $(M_n)$ is c.i.d. with
respect to the filtration
\[
\mathcal{G}_0=\{\varnothing,\Omega\}, \qquad\mathcal{G}_n=
\mathcal {F}_n\vee \sigma(R_{n+1})=\sigma(M_1,
\ldots,M_n,R_1,\ldots,R_n,R_{n+1}).
\]
We next prove this fact. The c.i.d. property has been recalled in
Section~\ref{kh8u5rf}.

%
%le1 #&#
\begin{lem}\label{c4p}
$(M_n)$ is c.i.d. with respect to $(\mathcal{G}_n)$ if and only if
%
%
%e1 #&#
\begin{equation}
\label{df7mp} \Lambda_{n+1}=\Lambda_n \biggl(1-
\frac{R_{n+1}-\beta}{c+\sum_{i=1}^{n+1}R_i} \biggr)\qquad \mbox{a.s. for all }n\geq0.
\end{equation}
In particular, $(M_n)$ is c.i.d. with respect to $(\mathcal{G}_n)$
if $\beta=0$ or if $R_n=1$ for all $n\geq1$. [In these cases, in
fact, condition \eqref{df7mp} is trivially true.]
\end{lem}

\begin{pf}
We just give a sketch of the proof. Suppose
%
%
%e2 #&#
\begin{equation}
\label{mag} M_{n+2}(B)\mid\mathcal{G}_n\sim
M_{n+1}(B)\mid\mathcal{G}_n \qquad\mbox{for each }n\geq0\mbox{ and
}B\in\mathcal{B}.
\end{equation}
Conditionally on $\mathcal{G}_n$, the r.m.'s $M_{n+1}$ and $M_{n+2}$
are both \emph{completely} r.m.'s. Hence, condition \eqref{mag} implies
\[
M_{n+2}\mid\mathcal{G}_n\sim M_{n+1}\mid
\mathcal{G}_n \qquad\mbox{for each }n\geq0.
\]
In turn, given $n\geq0$ and $A\in\Sigma$, the previous condition
yields
\begin{eqnarray*}
P (M_{n+3}\in A\mid\mathcal{G}_n )&=&E \bigl\{P
(M_{n+3}\in A\mid\mathcal{G}_{n+1} )\mid\mathcal{G}_n
\bigr\} \\
&=&E \bigl\{P (M_{n+2}\in A\mid\mathcal{G}_{n+1} )\mid
\mathcal{G}_n \bigr\}
\\
&=&P (M_{n+2}\in A\mid\mathcal{G}_n )=P (M_{n+1}\in
A\mid\mathcal{G}_n ) \qquad\mbox{a.s.}
\end{eqnarray*}
Hence, $M_{n+3}\mid\mathcal{G}_n\sim M_{n+1}\mid\mathcal{G}_n$ for
each $n\geq0$. Iterating this argument, one obtains
$M_k\mid\mathcal{G}_n\sim M_{n+1}\mid\mathcal{G}_n$ for all
$k>n\geq
0$. Therefore, condition \eqref{mag} is equivalent to $(M_n)$ being
c.i.d. with respect to $(\mathcal{G}_n)$. We next prove that
\eqref{df7mp} $\Leftrightarrow$ \eqref{mag}.

Fix $n\geq0$ and $B\in\mathcal{B}$. It can be assumed $m(B)>0$. Since
$R_{n+1}$ is independent of
$(M_1,\ldots,M_n,M_{n+1},R_1,\ldots,R_n)$, then
\[
P (M_{n+1}\in A\mid\mathcal{G}_n )=P (M_{n+1}\in A
\mid\mathcal{F}_n ) \qquad\mbox{a.s. for all }A\in\Sigma.
\]
Thus, for
each $t\in\mathbb{R}$,
\begin{eqnarray}
&&E \bigl\{e^{t M_{n+1}(B)}\mid\mathcal{G}_n \bigr\}\nonumber\\
&&\qquad=E \bigl
\{e^{t
M_{n+1}(B)}\mid\mathcal{F}_n \bigr\}
\nonumber\\
&&\qquad=\exp{ \bigl(m(B) \bigl(e^t-1 \bigr) \Lambda_n} \bigr)
\prod_{x\in S_n\cap
B} \biggl\{1+ \bigl(e^t-1 \bigr)
\frac{-\beta+\sum_{i=1}^nR_iM_i\{x\}}{c+\sum_{i=1}^nR_i} \biggr\}\nonumber\\
\eqntext{\mbox{a.s.},}
\end{eqnarray}
where the second equality is because $M_{n+1}\mid\mathcal{F}_n\sim
\operatorname{Be}P(\nu_n)$. Similarly,
\begin{eqnarray*}
&&E \bigl\{e^{t M_{n+2}(B)}\mid\mathcal{G}_n \bigr\}\nonumber\\
&&\qquad=E \bigl\{E
\bigl(e^{t
M_{n+2}(B)}\mid\mathcal{G}_{n+1} \bigr)\mid
\mathcal{G}_n \bigr\}
\nonumber\\
&&\qquad=\exp{ \bigl(m(B) \bigl(e^t-1 \bigr) \Lambda_{n+1}} \bigr)
\\
&&\qquad\quad{}\times E \biggl\{\prod_{x\in
S_{n+1}\cap B} \biggl(1+ \bigl(e^t-1
\bigr) \frac{-\beta+\sum_{i=1}^{n+1}R_iM_i\{
x\}
}{c+\sum_{i=1}^{n+1}R_i} \biggr)\Bigm|\mathcal{G}_n \biggr\}
\qquad\mbox{a.s.}
\end{eqnarray*}
Finally, after some computations, one obtains
\begin{eqnarray*}
&&E \biggl\{\prod_{x\in S_{n+1}\cap B} \biggl(1+
\bigl(e^t-1 \bigr) \frac{-\beta
+\sum_{i=1}^{n+1}R_iM_i\{x\}}{c+\sum_{i=1}^{n+1}R_i} \biggr)\Bigm|\mathcal
{G}_n \biggr\}
\\
&&\qquad=\exp{ \biggl(m(B) \bigl(e^t-1 \bigr) \Lambda_n
\frac{R_{n+1}-\beta
}{c+\sum_{i=1}^{n+1}R_i} \biggr)}\\
&&\qquad\quad{}\times \prod_{x\in S_n\cap B} \biggl\{1+
\bigl(e^t-1 \bigr) \frac{-\beta+\sum_{i=1}^nR_iM_i\{x\}}{c+\sum_{i=1}^nR_i} \biggr\}\qquad \mbox{a.s.}
\end{eqnarray*}
Thus, condition \eqref{df7mp} amounts to
$E \{e^{t M_{n+2}(B)}\mid\mathcal{G}_n \}=E \{e^{t
M_{n+1}(B)}\mid\mathcal{G}_n \}$
a.s. for each $t\in\mathbb{R}$, that is, conditions \eqref{df7mp}
and \eqref{mag} are equivalent.
\end{pf}

%s4 #&#
\section{Asymptotic behavior of $L_n$}\label{are}

Let $N_i$ be the number of new dishes tried by customer $i$, that is,
\[
N_i= \operatorname{card}(S_i\setminus
S_{i-1}) \qquad\mbox{with } S_0=\varnothing.
\]
Note that $N_i$ is $\mathcal{F}_i$-measurable and $N_i\mid\mathcal
{F}_{i-1}\sim\operatorname{Poi}(\Lambda_{i-1})$.

This section is devoted to
\[
L_n= \operatorname{card}(S_n)=\sum
_{i=1}^nN_i,
\]
the number of dishes experimented by the first $n$ customers. Our main
tool is the following technical lemma.

%
%le2 #&#
\begin{lem}\label{basic123l}
There is a function $h\dvtx(0,\infty)\rightarrow\mathbb{R}$ such that
\begin{eqnarray}
\sup_{x\geq c+u}\bigl|x h(x)\bigr|<\infty\quad\mbox{and} \quad\Lambda_n=
\alpha \frac{\Gamma(c+1)}{\Gamma(c+\beta)} \frac{1+h (c+\sum_{i=1}^nR_i )}{ (c+\sum_{i=1}^nR_i )^{1-\beta}}\nonumber \\
\eqntext{\mbox {for all }n\geq1.}
\end{eqnarray}
In particular,
%
%
%e3 #&#
\begin{equation}
\label{mc5r} \Lambda_n\leq\frac{D}{n^{1-\beta}}\quad \mbox{and}\quad |{\Lambda
_{n+1}-\Lambda_n }|\leq\frac{D}{n^{2-\beta}}\qquad \mbox{for all }n
\geq1,
\end{equation}
where $D$ is a suitable constant (nonrandom and not depending on $n$).
\end{lem}

\begin{pf}
Just note that $\frac{\Gamma(x+\beta)}{\Gamma(x+1)}=x^{\beta
-1}(1+h(x))$, with $h$ as required, for all $x>\max
{(0,-\beta
)}$; see, for example, formula (6.1.47) of \cite{AS}. To prove \eqref
{mc5r}, let $v=\min{(u,c+u)}$. Since $c+u>c+\beta>0$, then $v>0$.
Hence, \eqref{mc5r} follows from\looseness=-1
\[
c+\sum_{i=1}^nR_i\geq c+n
u=c+u+(n-1) u\geq n v.
\]\looseness=0
\upqed\end{pf}

Let $L=\sup_nL_n$ be the number of dishes tried by some customer. A
first consequence of Lemma~\ref{basic123l} is that $\beta<0$ implies
$L<\infty$ a.s.

%
%le3 #&#
\begin{lem}\label{pino78}
$P (N_i>\frac{1}{1-\beta}$ infinitely often$ )=0$. Moreover,
$E (e^L )<\infty$ if $\beta<0$.
\end{lem}

\begin{pf}
Fix an integer $k\geq1$. Since $N_{i+1}\mid\mathcal{F}_i\sim\operatorname{Poi}(\Lambda_i)$,
\[
P(N_{i+1}\geq k)=E \bigl\{P(N_{i+1}\geq k\mid
\mathcal{F}_i) \bigr\} =E \biggl\{e^{-\Lambda_i}\sum
_{j\geq k}\frac{\Lambda_i^j}{j!} \biggr\}\leq \frac
{E(\Lambda_i^k)}{k !}.
\]
By Lemma~\ref{basic123l}, $E(\Lambda_i^k)= \mathrm{O}(i^{-(1-\beta
)k})$. Let
\[
k=1+\max \bigl\{j\in\mathbb{Z}\dvtx j\leq1/(1-\beta) \bigr\}.
\]
Since $k (1-\beta)>1$, one obtains $\sum_iP (N_i> 1/(1-\beta
)
)=\sum_iP(N_i\geq k)<\infty$. Next, suppose $\beta<0$. By Lemma~\ref
{basic123l}, $\Lambda_n\leq D n^{\beta-1}$ for some constant $D$.
Letting $H=(e-1) D$ and noting that $E (e^{N_{n+1}}\mid\mathcal
{F}_n )=e^{\Lambda_n (e-1)}$ a.s., one obtains
\begin{eqnarray*}
E \bigl(e^{L_{n+1}} \bigr)&=&E \bigl\{e^{L_n} E
\bigl(e^{N_{n+1}}\mid \mathcal {F}_n \bigr) \bigr\}=E \bigl
\{e^{L_n} e^{\Lambda_n (e-1)} \bigr\}\leq E \bigl(e^{L_n} \bigr)
e^{H n^{\beta-1}}
\\
&\leq& E \bigl(e^{L_{n-1}} \bigr) e^{H (n-1)^{\beta-1}}e^{H n^{\beta
-1}}\leq\cdots
\leq E \bigl(e^{L_1} \bigr) e^{H \sum_{j=1}^nj^{\beta-1}}.
\end{eqnarray*}
Thus, $\beta<0$ and $E (e^{L_1} )=E (e^{N_1}
)<\infty$ yield
\[
E \bigl(e^L \bigr)=\sup_nE
\bigl(e^{L_n} \bigr)\leq E \bigl(e^{L_1} \bigr)
e^{H \sum_{j=1}^\infty j^{\beta-1}}<\infty.
\]
\upqed\end{pf}

In view of Lemma~\ref{pino78}, if $\beta<0$ there is a random index $N$
such that $L_n=L_N$ a.s. for all $n\geq N$. The situation is quite
different if $\beta\in[0,1)$. In this case, the a.s. behavior of
$L_n$ for large $n$ can be determined by a simple martingale argument.

In the rest of this section, we let $\beta\in[0,1)$. Define
\[
\overline{R}_n=\frac{1}{n}\sum_{i=1}^nR_i
\]
and suppose that
%
%
%e4 #&#
\begin{equation}
\label{bgbtt} \overline{R}_n\stackrel{\mathrm{a.s.}}\longrightarrow r\qquad
\mbox{for some constant }r.
\end{equation}
Since $R_i\geq u$ for all $i$, then $r\geq u>0$. Define also
\begin{eqnarray*}
\lambda(\beta)&=&\frac{\alpha c}{r} \qquad\mbox{if } \beta=0\quad \mbox {and}\quad \lambda(
\beta)=\frac{\alpha\Gamma(c+1)}{\Gamma(c+\beta)} \frac{1}{\beta r^{1-\beta}}\qquad \mbox{if } \beta\in(0,1),
\\
a_n(\beta)&=&\log{n} \qquad\mbox{if } \beta=0 \quad\mbox{and}\quad a_n(
\beta )=n^\beta\qquad\mbox{if } \beta\in(0,1).
\end{eqnarray*}

%
%th4 #&#
\begin{thmm}\label{thmln}
If $\beta\in[0,1)$ and condition \eqref{bgbtt} holds, then
\[
\frac{L_n}{a_n(\beta)}\stackrel{a.s.}\longrightarrow\lambda(\beta).
\]
\end{thmm}

\begin{pf}
By Lemma~\ref{basic123l}, $\Lambda_j=\alpha\frac{\Gamma
(c+1)}{\Gamma
(c+\beta)} (c+\sum_{i=1}^jR_i )^{\beta-1}  \{1+h
(c+\sum_{i=1}^jR_i ) \}$ where the function $h$ satisfies $
|{h(x)}|\leq(k/x)$ for all $x\geq c+u$ and some constant~$k$. Write
\begin{eqnarray}
\frac{\sum_{j=1}^{n-1}\Lambda_j}{a_n(\beta)} =\alpha\frac{\Gamma(c+1)}{\Gamma(c+\beta)} \frac{\sum_{j=1}^{n-1}j^{\beta-1}  ({c}/{j}+\overline{R}_j
)^{\beta
-1}}{a_n(\beta)}+D_n,
\nonumber\\
\eqntext{\displaystyle\mbox{where } D_n=\alpha\frac{\Gamma(c+1)}{\Gamma(c+\beta)} \frac{\sum_{j=1}^{n-1} (c+\sum_{i=1}^jR_i )^{\beta
-1}h
(c+\sum_{i=1}^jR_i )}{a_n(\beta)}.}
\end{eqnarray}
In view of \eqref{bgbtt}, one obtains $D_n\stackrel
{\mathrm{a.s.}}\longrightarrow
0$ and $\frac{\sum_{j=1}^{n-1}\Lambda_j}{a_n(\beta)}\stackrel
{\mathrm{a.s.}}\longrightarrow\lambda(\beta)$. Next, define
\[
T_0=0\quad \mbox{and}\quad T_n=\sum_{j=1}^n
\frac{N_j-E(N_j\mid\mathcal
{F}_{j-1})}{a_j(\beta)}=\sum_{j=1}^n
\frac{N_j-\Lambda
_{j-1}}{a_j(\beta)}.
\]
Then, $(T_n)$ is a martingale with respect to $(\mathcal{F}_n)$ and
\begin{eqnarray*}
E \bigl(T_n^2 \bigr)&=&\sum_{j=1}^n
\frac{E \{(N_j-\Lambda_{j-1})^2 \}
}{a_j(\beta)^2}=\sum_{j=1}^n
\frac{E \{E ((N_j-\Lambda
_{j-1})^2\mid\mathcal{F}_{j-1} ) \}}{a_j(\beta)^2}\\
&=&\sum_{j=1}^n
\frac{E(\Lambda_{j-1})}{a_j(\beta)^2}.
\end{eqnarray*}
Since $E(\Lambda_j)=\mathrm{O}(j^{-(1-\beta)})$, then $\sup_nE(T_n^2)=\sum_{j=1}^\infty\frac{E(\Lambda_{j-1})}{a_j(\beta)^2}<\infty$. Thus, $T_n$
converges a.s., and Kronecker's lemma implies
\begin{eqnarray*}
\lim_n\frac{L_n}{a_n(\beta)}&=&\lim_n
\frac{\sum_{j=1}^nN_j}{a_n(\beta
)}=\lim_n\frac{\sum_{j=1}^n\Lambda_{j-1}}{a_n(\beta)}\\
&=&\lim
_n\frac
{\Lambda_0+\sum_{j=1}^{n-1}\Lambda_j}{a_n(\beta)}=\lambda(\beta)
\qquad\mbox{a.s.}
\end{eqnarray*}
\upqed\end{pf}

In view of Theorem~\ref{thmln}, as far as $\beta\in[0,1)$ and the
weights $R_n$ meet the
SLLN, $L_n$ essentially behaves for large $n$ as in the standard IBP
model. The only difference is that the limit constant
$\lambda(\beta)$ depends on $r$ as well. (In the standard IBP one has
$r=1$.) Note also that, the $R_n$ being independent, a
sufficient condition for \eqref{bgbtt} is
\[
\sup_nE \bigl(R_n^2 \bigr)<\infty\quad
\mbox{and}\quad \frac{\sum_{i=1}^nE(R_i)}{n}\longrightarrow r.
\]

We next turn to the limiting distribution of $L_n$. To get something,
stronger conditions on the $R_n$ are to be requested.

%
%th5 #&#
\begin{thmm}\label{abg}
If $\beta\in[0,1)$ and
%
%
%e5 #&#
\begin{equation}
\label{jkh6t} \overline{R}_n\stackrel{a.s.}\longrightarrow r
\quad\mbox{and}\quad \frac
{\sum_{j=1}^nj^{\beta-1} E{|\overline{R}_j-r}|}{\sqrt{a_n(\beta
)}}\longrightarrow0
\end{equation}
for some constant $r$, then
\[
\sqrt{a_n(\beta)} \biggl\{\frac{L_n}{a_n(\beta)}-\lambda(\beta ) \biggr
\} \longrightarrow\mathcal{N} \bigl(0, \lambda(\beta) \bigr)\qquad \mbox{stably}.
\]
\end{thmm}

\begin{pf}
We first prove that
%
%
%e6 #&#
\begin{equation}
\label{uje} \sqrt{a_n(\beta)} \biggl\{\frac{\sum_{j=1}^n\Lambda
_{j-1}}{a_n(\beta
)}-\lambda(
\beta) \biggr\}\stackrel{P}\longrightarrow0.
\end{equation}
By Lemma~\ref{basic123l} and some calculations, condition \eqref{uje}
is equivalent to
\[
Y_n:=\frac{\sum_{j=1}^{n-1}  \{ (c+\sum_{i=1}^jR_i
)^{\beta
-1}-(r j)^{\beta-1} \}}{\sqrt{a_n(\beta)}}\stackrel {P}\longrightarrow0.
\]
Let $v=\min{(u,c+u)}$. Then, $v>0$, $r\geq u\geq v$ and $c+\sum_{i=1}^jR_i\geq v j$; see the proof of Lemma~\ref{basic123l}. Hence,
one can estimates as follows:
\begin{eqnarray*}
E\Biggl|(r j)^{\beta-1}- \Biggl(c+\sum_{i=1}^jR_i
\Biggr)^{\beta-1}\Biggr|&\leq& \frac
{E| (c+\sum_{i=1}^jR_i )^{1-\beta}-(r j)^{1-\beta}|}{(v
j)^{2(1-\beta)}}
\\
&\leq&\frac{1}{(v j)^{2(1-\beta)}} \frac{1-\beta}{(v j)^{\beta}} E\Biggl|c+\sum_{i=1}^jR_i-r
j\Biggr| \\
&\leq&\frac{1-\beta}{v^{2-\beta}} \biggl\{\frac{|c|}{j^{2-\beta
}}+\frac{E|\overline{R}_j-r|}{j^{1-\beta}} \biggr\}.
\end{eqnarray*}
Thus, condition \eqref{jkh6t} implies $E|Y_n|\rightarrow0$. This
proves condition \eqref{uje}.

Next, define
\[
U_n=\sqrt{a_n(\beta)} \biggl\{\frac{L_n}{a_n(\beta)}-
\frac{\sum_{j=1}^n\Lambda_{j-1}}{a_n(\beta)} \biggr\}=\frac{\sum_{j=1}^n(N_j-\Lambda
_{j-1})}{\sqrt{a_n(\beta)}}.
\]
In view of \eqref{uje}, it suffices to show that $U_n\longrightarrow
\mathcal{N} (0, \lambda(\beta) )$ stably. To this end, for
$n\geq1$ and $j=1,\ldots,n$, define
\[
U_{n,j}=\frac{N_j-\Lambda_{j-1}}{\sqrt{a_n(\beta)}}, \qquad\mathcal {R}_{n,0}=
\mathcal{F}_0\quad \mbox{and}\quad \mathcal{R}_{n,j}=
\mathcal{F}_j.
\]
Then, $E (U_{n,j}\mid\mathcal{R}_{n,j-1})=0$ a.s., $\mathcal
{R}_{n,j}\subset\mathcal{R}_{n+1,j}$ and $U_n=\sum_jU_{n,j}$. Thus, by
the martingale CLT, $U_n\longrightarrow\mathcal{N} (0, \lambda
(\beta) )$ stably provided
\begin{eqnarray*}
&&\mbox{(i)}\quad\sum_{j=1}^nU_{n,j}^2
\stackrel{P}\longrightarrow\lambda (\beta ),\qquad \mbox{(ii)}\quad\max
_{1\leq j\leq n}|U_{n,j}| \stackrel {P}\longrightarrow0,\\
&&\mbox{(iii)}\quad\sup_nE \Bigl\{ \max_{1\leq j\leq
n}U_{n,j}^2
\Bigr\}<\infty;
\end{eqnarray*}
see, for example, Theorem~3.2, page 58, of \cite{HH}. Let
\[
H_j=(N_j-\Lambda_{j-1})^2
\quad\mbox{and}\quad D_n=\frac{\sum_{j=1}^n \{H_j-E(H_j\mid\mathcal{F}_{j-1}) \}}{a_n(\beta)}.
\]
By Kronecker's lemma and the same martingale argument used in the proof
of Theorem~\ref{thmln}, $D_n\stackrel{\mathrm{a.s.}}\longrightarrow0$. Since
$\overline{R}_j\stackrel{\mathrm{a.s.}}\longrightarrow r$ and $E(H_j\mid
\mathcal
{F}_{j-1})=\Lambda_{j-1}$ a.s., then
\[
\sum_{j=1}^nU_{n,j}^2=
\frac{\sum_{j=1}^nH_j}{a_n(\beta)}=D_n+\frac
{\sum_{j=1}^n\Lambda_{j-1}}{a_n(\beta)}\stackrel{\mathrm{a.s.}}
\longrightarrow \lambda (\beta).
\]
This proves condition (i). As to (ii), fix $k\geq1$ and note that
\begin{eqnarray}
\max_{1\leq j\leq n}U_{n,j}^2\leq
\frac{\max_{1\leq j\leq
k}H_j}{a_n(\beta)}+\max_{k< j\leq n} \frac{H_j}{a_j(\beta)}\leq
\frac
{\max_{1\leq j\leq k}H_j}{a_n(\beta)}+\sup_{j>k} \frac
{H_j}{a_j(\beta
)}\nonumber\\
 \eqntext{\mbox{for }n>k.}
\end{eqnarray}
Hence, $\limsup_n\max_{1\leq j\leq n}U_{n,j}^2\leq\limsup_n\frac
{H_n}{a_n(\beta)}$ and condition (ii) follows from
\[
\frac{H_n}{a_n(\beta)}=\frac{\sum_{j=1}^nH_j}{a_n(\beta)}-\frac
{\sum_{j=1}^{n-1}H_j}{a_n(\beta)}\stackrel{\mathrm{a.s.}}
\longrightarrow0.
\]
Finally, condition (iii) is an immediate consequence of Lemma~\ref
{basic123l} and
\[
E \Bigl\{\max_{1\leq j\leq n}U_{n,j}^2 \Bigr\}\leq
\frac{\sum_{j=1}^nE(H_j)}{a_n(\beta)}=\frac{\sum_{j=1}^nE(\Lambda
_{j-1})}{a_n(\beta)}.
\]
\upqed\end{pf}

Note that, letting $R_n=1$ for all $n$, Theorem~\ref{abg} provides the
limiting distribution of $L_n$ in the standard IBP model.

For Theorem~\ref{abg} to apply, condition \eqref{jkh6t} is to be
checked. We now give conditions for \eqref{jkh6t}. In particular,
\eqref
{jkh6t} is automatically true whenever $\sup_nE(R_n^2)<\infty$ and
$E(R_n)=r$ for all $n$.

%
%le6 #&#
\begin{lem}
Condition \eqref{jkh6t} holds provided $\beta\in[0,1)$ and
%
%
%e7 #&#
\begin{equation}
\label{spib56t} \sup_nE \bigl(R_n^2
\bigr)<\infty\quad\mbox{and}\quad \sqrt{n^{\beta} \log{n}} \bigl\{E(
\overline{R}_n)-r \bigr\}\longrightarrow0.
\end{equation}
\end{lem}

\begin{pf}
Let $a=\sup_nE(R_n^2)$. Because of \eqref{spib56t}, $E(\overline
{R}_n)\rightarrow r$. Thus, $\overline{R}_n\stackrel
{\mathrm{a.s.}}\longrightarrow r$ since $a<\infty$ and $(R_n)$ is independent. Moreover,
\begin{eqnarray*}
E|\overline{R}_j-r|&\leq& E\bigl|\overline{R}_j-E(
\overline{R}_j)\bigr|+ \bigl|E(\overline{R}_j)-r\bigr|\leq\sqrt{
\operatorname{var}(\overline{R}_j)}+ \bigl|E(\overline{R}_j)-r\bigr|
\\
&\leq&\sqrt{a/j}+\bigl|E(\overline{R}_j)-r\bigr|.
\end{eqnarray*}
Hence, the second part of condition \eqref{jkh6t} follows from the
above inequality and condition \eqref{spib56t}.
\end{pf}

A last remark is in order. Fix a set $B\in\mathcal{B}$ and define
\[
L_n(B)= \operatorname{card}(B\cap S_n)
\]
to be the number of dishes, belonging to $B$, tried by the first $n$
customers. The same arguments used for $L_n=L_n(\mathcal{X})$ apply to
$L_n(B)$ and allow us to extend Theorems~\ref{thmln}--\ref{abg} as follows.

%
%th7 #&#
\begin{thmm}
Let $\beta\in[0,1)$ and $B\in\mathcal{B}$. If condition \eqref{bgbtt}
holds, then
\[
\frac{L_n(B)}{a_n(\beta)}\stackrel{a.s.}\longrightarrow m(B)\lambda (\beta).
\]
Moreover, under condition \eqref{jkh6t}, one obtains
\[
\sqrt{a_n(\beta)} \biggl\{\frac{L_n(B)}{a_n(\beta)}-m(B)\lambda (\beta )
\biggr\}\longrightarrow\mathcal{N} \bigl(0, m(B)\lambda(\beta ) \bigr)\qquad
\mbox{stably}.
\]
\end{thmm}

\begin{pf}
Let $N_i(B)$ denote the number of new dishes, belonging to $B$, tried
by customer $i$. Then, $L_n(B)=\sum_{i=1}^nN_i(B)$ and $N_{i+1}(B)\mid
\mathcal{F}_i\sim\break \operatorname{Poi}(m(B) \Lambda_i)$. Therefore, it suffices to
repeat the proofs of Theorems \ref{thmln}--\ref{abg} with $N_i(B)$ in
the place of $N_i$ and $m(B) \Lambda_i$ in the place of $\Lambda_i$.
\end{pf}

%s5 #&#
\section{Asymptotic behavior of $\overline{K}_n$}

%s5.1 #&#
\subsection{The result}\label{v4gr} Let $K_i=M_i(\mathcal{X})$ be the
number of dishes experimented by customer $i$ and
\[
\overline{K}_n=\frac{1}{n}\sum_{i=1}^nK_i
\]
the mean number of dishes tried by each of the first $n$ customers.

In IBP-type models, $\overline{K}_n$ is a meaningful quantity. One
reason is the following. If the parameters $m, \alpha, \beta$ and $c$
are unknown, $E (K_{n+1}\mid\mathcal{F}_n )$ cannot be
evaluated in closed form. Then, $\overline{K}_n$ could be used as an
empirical predictor for the next random variable $K_{n+1}$. Such
prediction is consistent whenever
\[
V_n:=\overline{K}_n-E (K_{n+1}\mid
\mathcal{F}_n )\stackrel {P}\longrightarrow 0.
\]
But this is usually true. For instance,
$V_n\stackrel{\mathrm{a.s.}}\longrightarrow0$ if the sequence $(K_n)$ is
c.i.d. with respect to $(\mathcal{F}_n)$; see \cite{BPR04} and
\cite{BCPR09}. In general, the higher the convergence rate of $V_n$,
the better $\overline{K}_n$ as a predictor of $K_{n+1}$.

Under some conditions, $\overline{K}_n\stackrel{\mathrm{a.s.}}\longrightarrow
Z$ for some real random variable $Z$. Thus, two random centerings
for $\overline{K}_n$ should be considered. One (and more natural) is
$Z$, while the other is $E (K_{n+1}\mid\mathcal{F}_n )$, to
evaluate the performances of $\overline{K}_n$ as a predictor of
$K_{n+1}$. Taking $\sqrt{n}$ as a norming factor, this leads to
investigate
\[
\sqrt{n} \{\overline{K}_n-Z \} \quad\mbox{and}\quad \sqrt{n} V_n.
\]
The limiting distributions of these quantities are
provided by the next result.

%
%th8 #&#
\begin{thmm}\label{cert5r}
Suppose $\beta<1/2$ and
\[
\sup_nR_n\leq b,\qquad E(R_n)
\longrightarrow r,\qquad E \bigl(R_n^2 \bigr)\longrightarrow q,
\]
for some constants $b, r, q$. Then
\[
\overline{K}_n\stackrel{a.s.}\longrightarrow Z \quad\mbox{and}\quad
\frac{1}{n} \sum_{i=1}^nK_i^2
\stackrel {a.s.}\longrightarrow Q,
\]
where $Z$ and $Q$ are real random variables such that $Z^2<Q$ a.s.
Moreover,
\begin{eqnarray}
\sqrt{n} \{\overline{K}_n-Z \}&\longrightarrow&\mathcal {N} \bigl(0,
\sigma^2 \bigr) \qquad\mbox{stably and}\nonumber
\\
\sqrt{n} \bigl\{\overline{K}_n-E (K_{n+1}\mid\mathcal
{F}_n ) \bigr\}&\longrightarrow&\mathcal{N} \bigl(0, \tau^2
\bigr) \qquad\mbox{stably,}\nonumber
\\
\eqntext{\mbox{where } \displaystyle\sigma^2=\frac{2q-r^2}{r^2} \bigl(Q-Z^2
\bigr), \tau ^2=\frac{q-r^2}{r^2} \bigl(Q-Z^2 \bigr).}
\end{eqnarray}
If $R_n=1$ for all $n$, the previous results hold for $\beta<1$ (and
not only for $\beta<1/2$).
\end{thmm}

Theorem~\ref{cert5r} is a consequence of Theorem~1 of \cite{BCPR11}.
The proof, even if conceptually simple, is technically rather hard.

Theorem~\ref{cert5r} fails, as it stands, for $\beta\in[1/2, 1)$. Let
$\mu_n$ denote the probability distribution of the random variable
$\sqrt{n}  \{\overline{K}_n-Z \}$. The sequence $(\mu_n)$
might be not tight if $\beta\in(1/2, 1)$. For instance, $(\mu_n)$ is
not tight if $\beta\in(1/2, 1)$ and $R_n=r$ for all
$n$, where $r$ is any constant such that $r\neq1$. If $\beta=1/2$,
instead, $(\mu_n)$ is tight, but the possible
limit laws are not mixtures of centered Gaussian distributions. Thus,
even if $\sqrt{n}  \{\overline{K}_n-Z \}$ converges stably,
the limit kernel is not $\mathcal{N} (0, \sigma^2 )$.

Since $q\geq r^2$ and $Q>Z^2$ a.s., then $\sigma^2>0$ a.s. Hence,
$\mathcal{N} (0, \sigma^2 )$ is a nondegenerate kernel.
Instead, $\mathcal{N} (0, \tau^2 )$ may be degenerate. In
fact, if $q=r^2$ then
$\mathcal{N} (0, \tau^2 )=\mathcal{N} (0,0)=\delta_0$.
Thus, for $q=r^2$, Theorem~\ref{cert5r} yields
$\sqrt{n} V_n\stackrel{P}\longrightarrow0$.

The convergence rate of $V_n$ is $n^{-1/2}$ when $q>r^2$. Such a
rate is even higher if $q=r^2$, since $\sqrt{n} V_n\stackrel
{P}\longrightarrow0$. Overall,
$\overline{K}_n$ seems to be a good predictor of $K_{n+1}$ for large
$n$.

Among other things, Theorem~\ref{cert5r} can be useful to get
asymptotic confidence bounds for $Z$. Define in fact
\[
\widehat{\sigma}_n^2= \biggl\{\frac{(2/n) \sum_{i=1}^nR_i^2}{\overline
{R}_n^2}-1
\biggr\} \Biggl\{\frac{1}{n} \sum_{i=1}^nK_i^2-
\overline {K}_n^2 \Biggr\}.
\]
Since $\widehat{\sigma}_n^2\stackrel{\mathrm{a.s.}}\longrightarrow\sigma^2$ and
$\sigma^2>0$ a.s., one obtains
\[
I_{\{\widehat{\sigma}_n>0\}}\frac{\sqrt{n}  \{\overline
{K}_n-Z
\}}{\widehat{\sigma}_n}\longrightarrow\mathcal{N}(0, 1)\qquad
\mbox{stably}.
\]
Thus, $\overline{K}_n\pm\frac{u_a}{\sqrt{n}} \widehat{\sigma}_n$
provides an asymptotic confidence interval for $Z$ with
(approximate) level $1-a$, where $u_a$ is such that
$\mathcal{N}(0,1)(u_a, +\infty)=a/2$.

Theorem~\ref{cert5r} works if $\beta\in[0,1)$ and $R_n=1$ for all $n$,
that is, it applies to the standard IBP model. Also, in this case, the
convergence
rate of $V_n$ is greater than $n^{-1/2}$ (since $q=r^2=1$). Hence,
$\overline{K}_n$ is
a good (asymptotic) predictor of $K_{n+1}$.

%s5.2 #&#
\subsection{The proof} We begin with a couple of results from \cite
{BCPR11}. Let $(X_n)$ be a sequence of real integrable random variables,
adapted to a filtration $(\mathcal{U}_n)$, and let
\[
\overline{X}_n=\frac{1}{n}\sum_{i=1}^nX_i\quad
\mbox{and}\quad Z_n=E (X_{n+1}\mid\mathcal{U}_n ).
\]

%
%le9 #&#
\begin{lem}\label{iam7t}
If $\sum_nn^{-2}E(X_n^2)<\infty$ and $Z_n\stackrel
{a.s.}\longrightarrow
Z$, for some real random variable $Z$, then
\[
\overline{X}_n\stackrel{a.s.}\longrightarrow Z \quad\mbox{and}\quad n\sum
_{k\geq n}\frac{X_k}{k^2}\stackrel{a.s.}\longrightarrow
Z.
\]
\end{lem}

\begin{pf}
This is exactly Lemma~2 of \cite{BCPR11}.
\end{pf}

%
%th10 #&#
\begin{thmm}\label{j2a87p}
Suppose $(X_n^2)$ is uniformly integrable and

\begin{longlist}[(j)]
\item[(j)] $n^3E \{ (E(Z_{n+1}\mid\mathcal
{U}_n)-Z_n
)^2 \}\longrightarrow
0$.
\end{longlist}

Then $Z_n\stackrel{a.s.}\longrightarrow Z$ and $\overline
{X}_n\stackrel{a.s.}\longrightarrow Z$ for some real random variable
$Z$. Moreover,
\begin{eqnarray*}
\sqrt{n} \{\overline{X}_n-Z_n \}&\longrightarrow&\mathcal
{N} (0, U ) \qquad\mbox{stably}\mbox{ and}
\\
\sqrt{n} \{\overline{X}_n-Z \}&\longrightarrow&\mathcal {N} (0, U+V )\qquad
\mbox{stably}
\end{eqnarray*}
for some real random variables $U$ and $V$, provided

\begin{longlist}[(jjj)]
\item[(jj)] $E \{\sup_{k\geq1}\sqrt{k} |Z_{k-1}-Z_k|
\}<\infty$,

\item[(jjj)] $\frac{1}{n}\sum_{k=1}^n \{
X_k-Z_{k-1}+k(Z_{k-1}-Z_k) \}^2\stackrel{P}\longrightarrow
U$,

\item[(jv)] $n\sum_{k\geq n}(Z_{k-1}-Z_k)^2\stackrel
{a.s.}\longrightarrow
V$.
\end{longlist}
\end{thmm}

\begin{pf}
First note that $(Z_n)$ is a quasi-martingale because of (j) and it is
uniformly integrable for $(X_n^2)$ is uniformly integrable. Hence,
$Z_n\stackrel{\mathrm{a.s.}}\longrightarrow Z$. By Lemma~\ref{iam7t}, one also
obtains $\overline{X}_n\stackrel{\mathrm{a.s.}}\longrightarrow Z$. Next, assume
conditions (jj)--(jjj)--(jv). By Theorem~1 of \cite{BCPR11} (and the
subsequent remarks) it is enough to show that
\[
\sqrt{n} E \Bigl\{\sup_{k\geq n}|Z_{k-1}-Z_k|
\Bigr\} \longrightarrow 0\quad \mbox{and}\quad \frac{1}{\sqrt{n}} E \Bigl\{\max
_{1\leq k\leq
n}k |Z_{k-1}-Z_k| \Bigr\}
\longrightarrow0.
\]
Let $D_k=|Z_{k-1}-Z_k|$. Because of (jv),
\[
n D_n^2=n\sum_{k\geq n}D_k^2-
\frac{n}{n+1} (n+1)\sum_{k\geq
n+1}D_k^2
\stackrel{\mathrm{a.s.}}\longrightarrow 0.
\]
Thus $\sup_{k\geq n}\sqrt{k} D_k\stackrel{\mathrm{a.s.}}\longrightarrow
0$, and condition (jj) implies
\[
\sqrt{n} E \Bigl\{\sup_{k\geq n}D_k \Bigr\}\leq E
\Bigl\{\sup_{k\geq
n}\sqrt{k} D_k \Bigr\}
\longrightarrow 0.
\]
Further, for $1\leq h < n$, one obtains
\begin{eqnarray*}
E \Bigl\{\max_{1\leq k\leq
n}k D_k \Bigr\}&\leq& E \Bigl\{
\max_{1\leq k\leq h}k D_k \Bigr\}+\sqrt{n} E \Bigl\{\max
_{h<k\leq
n}\sqrt{k} D_k \Bigr\}\\
&\leq & E \Bigl\{\max
_{1\leq k\leq
h}k D_k \Bigr\}+\sqrt{n} E \Bigl\{\sup
_{k>h}\sqrt{k} D_k \Bigr\}.
\end{eqnarray*}
Hence, it suffices to note that
\begin{eqnarray*}
\limsup_n\frac{1}{\sqrt{n}} E \Bigl\{\max
_{1\leq k\leq
n}k D_k \Bigr\}&\leq& E \Bigl\{\sup
_{k>h}\sqrt{k} D_k \Bigr\} \qquad\mbox{for all }h\quad\mbox{and }\\
\lim_hE \Bigl\{\sup_{k>h}\sqrt{k}
D_k \Bigr\}&=&0.
\end{eqnarray*}
\upqed\end{pf}

Note that condition (j) is automatically true in case $(X_n)$ is c.i.d.
with respect to the filtration $(\mathcal{U}_n)$. We are now able to
prove Theorem~\ref{cert5r}.

\begin{pf*}{Proof of Theorem~\ref{cert5r}}
We apply Theorem~\ref
{j2a87p} with $X_n=K_n$ and $\mathcal{U}_n=\mathcal{F}_n$. Let
\[
J_n(x)=\frac{\sum_{i=1}^nR_iM_i\{x\}-\beta}{\sum_{i=1}^nR_i+c}\qquad \mbox {for }x\in\mathcal{X}.
\]
Note that
\[
\sum_{x\in S_n}J_n(x)=\frac{\sum_{i=1}^nR_i\sum_{x\in S_n}M_i\{x\}
-\beta
L_n}{\sum_{i=1}^nR_i+c} =
\frac{\sum_{i=1}^nR_iK_i-\beta L_n}{\sum_{i=1}^nR_i+c},
\]
and recall the notation
\[
\mathcal{G}_n=\mathcal{F}_n\vee\sigma(R_{n+1})=
\sigma(M_1,\ldots,M_n,R_1,
\ldots,R_n,R_{n+1}).
\]

\textit{Uniform integrability of} $(K_n^2)$. It suffices to show that
$\sup_nE \{e^{t K_n} \}<\infty$ for some $t>0$. In particular,
$(K_n^2)$ is uniformly integrable for $\beta<0$, since Lemma~\ref
{pino78} yields
\[
\sup_nE \bigl\{e^{K_n} \bigr\}\leq\sup
_nE \bigl\{e^{L_n} \bigr\} =E \bigl\{
e^L \bigr\}<\infty\qquad\mbox{if }\beta<0.
\]
Suppose $\beta\in[0, 1/2)$. Define $g(t)=e^t-1$ and
\[
W_n=\frac{\sum_{i=1}^nR_iK_i}{\sum_{i=1}^nR_i+c}.
\]
Arguing as in Lemma~\ref{c4p} and since $\Lambda_n\leq D n^{\beta-1}$
for some constant $D$, one obtains
\begin{eqnarray*}
E \bigl\{e^{t K_{n+1}}\mid\mathcal{G}_n \bigr\}&=&e^{g(t) \Lambda_n}
\prod_{x\in S_n} \bigl\{1+g(t) J_n(x) \bigr\}
\\
&\leq&\exp{ \biggl\{g(t) \Lambda_n+g(t)\sum
_{x\in S_n}J_n(x) \biggr\}}
\\
&\leq&\exp{ \biggl\{\frac{D g(t)}{n^{1-\beta}}+g(t) \frac{\sum_{i=1}^nR_iK_i-\beta L_n}{\sum_{i=1}^nR_i+c} \biggr\}}\\
&\leq&\exp{
\biggl\{ \frac{D g(t)}{n^{1-\beta}}+g(t) W_n \biggr\}} \qquad\mbox{a.s.}
\end{eqnarray*}
Hence, it is enough to show that $\sup_nE \{e^{t W_n} \}
<\infty
$ for some $t>0$. We first prove $E \{e^{t W_n} \}<\infty$ for
all $n\geq1$ and $t>0$, and subsequently $\sup_nE \{e^{t
W_n}
\}<\infty$ for a suitable $t>0$. Define $U_n=\frac{R_{n+1}}{\sum_{i=1}^{n+1}R_i+c}$. Since $U_n$ is $\mathcal{G}_n$-measurable,
\begin{eqnarray*}
E \bigl(e^{t W_{n+1}} \bigr)&=&E \bigl\{\exp{ \bigl(t W_n
(1-U_n) \bigr)} E \bigl(e^{t U_n K_{n+1}}\mid\mathcal{G}_n
\bigr) \bigr\}
\\
&\leq& E \biggl\{\exp{ \biggl(\frac{D g(t U_n)}{n^{1-\beta}} \biggr)} \exp { \bigl(t
W_n + \bigl(g(t U_n)-t U_n \bigr)
W_n \bigr)} \biggr\}.
\end{eqnarray*}
On noting that $U_n\leq b/(n u)$,
\[
E \bigl(e^{t W_{n+1}} \bigr)\leq\exp{ \biggl(\frac{D g (
{tb}/{(nu)} )}{n^{1-\beta}} \biggr)} E
\bigl\{e^{\{t+g({tb}/{(nu)})\}
W_n} \bigr\}.
\]
Iterating this procedure, one obtains
\[
E \bigl(e^{t W_{n+1}} \bigr)\leq a_n(t) E \bigl(e^{b_n(t) W_1}
\bigr)\qquad \mbox{for suitable constants }a_n(t)\mbox{ and
}b_n(t).
\]
Since $K_1\sim\operatorname{Poi}(\alpha)$ and $W_1=\frac{R_1}{R_1+c} K_1\leq
\frac{b}{u+c} K_1$, then $E (e^{b_n(t) W_1} )<\infty$. Hence,
$E \{e^{t W_n} \}<\infty$ for all $n\geq1$ and $t>0$. Observe
now that $g(z)\leq2 z$ and $g(z)-z\leq z^2$ for $z\in[0, 1/2]$.
Since $U_n\leq b/(n u)$, then $t U_n\leq1/2$ for $n\geq(2 b
t)/u$. Hence, if $t\in(0,1]$ and $n\geq(2 b)/u$, then
%
%
%e8 #&#
\begin{eqnarray}
\label{mb} E \bigl(e^{t W_{n+1}} \bigr)&\leq&\exp{ \biggl\{
\frac{2 D (b/u)
t}{n^{2-\beta}}} \biggr\} E \bigl\{\exp{ \bigl(t W_n+(t
U_n)^2W_n \bigr)} \bigr\}
\nonumber
\\[-8pt]
\\[-8pt]
\nonumber
&\leq&\exp{ \biggl\{\frac{D^* t}{n^{2-\beta}}} \biggr\} E \biggl\{ \exp{ \biggl(t
W_n \biggl(1+\frac{D^*}{n^2} \biggr) \biggr)} \biggr\},
\end{eqnarray}
where $D^*=\max{\{2 D (b/u), (b/u)^2\}}$. Take $t$ and $n_0$ such that
\[
t\in(0, 1/2], \qquad n_0\geq\frac{2 b}{u},\qquad \prod
_{j\geq n_0} \biggl(1+\frac{D^*}{j^2} \biggr)\leq2.
\]
Iterating inequality \eqref{mb}, one finally obtains
\[
E \bigl(e^{t W_{n+1}} \bigr)\leq\exp{ \biggl\{\sum
_{j\geq n_0}\frac{2
D^* t}{j^{2-\beta}} \biggr\}} E \bigl(e^{2 t W_{n_0}}
\bigr) \qquad\mbox {for each }n\geq n_0.
\]
Therefore $\sup_nE \{e^{t W_n} \}<\infty$, so that $(K_n^2)$
is uniformly integrable.

We now turn to condition (j). Since $M_{n+1}\mid\mathcal{F}_n\sim
\operatorname{Be}P(\nu_n)$,
\[
Z_n=E(K_{n+1}\mid\mathcal{F}_n)=
\nu_n(\mathcal{X}) =\Lambda_n+\sum
_{x\in S_n}J_n(x)=\Lambda_n+
\frac{\sum_{i=1}^nR_iK_i-\beta L_n}{\sum_{i=1}^nR_i+c}.
\]
On noting that $L_n=L_{n-1}+N_n$, a simple calculation yields
\[
Z_n-Z_{n-1}=\Lambda_n-\Lambda_{n-1}+
\frac
{R_n(K_n-Z_{n-1})+R_n\Lambda
_{n-1}-\beta N_n}{\sum_{i=1}^nR_i+c}.
\]

\textit{Condition} (j). Since $R_{n+1}$ is independent of $(M_1,\ldots,M_n,M_{n+1},R_1,\ldots,R_n)$,
\begin{eqnarray*}
E(K_{n+1}\mid\mathcal{G}_n)&=&E(K_{n+1}\mid
\mathcal{F}_n)=Z_n \quad\mbox {and}\\
 E(N_{n+1}\mid
\mathcal{G}_n)&=&E(N_{n+1}\mid\mathcal {F}_n)=
\Lambda_n \qquad\mbox{a.s.}
\end{eqnarray*}
It follows that
\begin{eqnarray*}
E (Z_{n+1}-Z_n\mid\mathcal{F}_n )&=&E \bigl\{E
(Z_{n+1}-Z_n\mid \mathcal{G}_n )\mid
\mathcal{F}_n \bigr\}\\
&=&E \biggl\{\Lambda _{n+1}-
\Lambda_n+\frac{(R_{n+1}-\beta) \Lambda_n}{\sum_{i=1}^{n+1}R_i+c}\Bigm|\mathcal{F}_n \biggr\}\qquad
\mbox{a.s.}
\end{eqnarray*}
Hence,
\begin{eqnarray*}
E \bigl\{ \bigl(E(Z_{n+1}\mid\mathcal{F}_n)-Z_n
\bigr)^2 \bigr\} &=&E \bigl\{ E (Z_{n+1}-Z_n\mid
\mathcal{F}_n )^2 \bigr\}
\\
&\leq&2 E \bigl\{(\Lambda_{n+1}-\Lambda_n)^2
\bigr\}+\frac{2 (b+
|\beta|)^2}{u^2} \frac{E(\Lambda_n^2)}{n^2}.
\end{eqnarray*}
By Lemma~\ref{basic123l}, $E \{(\Lambda_{n+1}-\Lambda_n)^2
\}
=\mathrm{O}(n^{2\beta-4})$ and $E(\Lambda_n^2)=\mathrm{O}(n^{2\beta-2})$.
Hence, condition (j) follows from $\beta<1/2$ (or equivalently
$4-2\beta>3$).

Having proved condition (j) and $(K_n^2)$ uniformly integrable, Theorem~\ref{j2a87p} yields $Z_n\stackrel{\mathrm{a.s.}}\longrightarrow Z$ and
$\overline
{K}_n\stackrel{\mathrm{a.s.}}\longrightarrow Z$ for some $Z$. We next prove
$(1/n) \sum_{i=1}^nK_i^2\stackrel{\mathrm{a.s.}}\longrightarrow Q$ for some $Q$
such that $Q>Z^2$ a.s. Recall that
\[
\sup_nE \bigl(K_n^4 \bigr)\leq
\frac{4!}{t^4} \sup_nE \bigl(e^{t K_n} \bigr)<
\infty \qquad\mbox{for a suitable }t>0.
\]
Hence, by Lemma~\ref{iam7t}, $(1/n) \sum_{i=1}^nK_i^2\stackrel
{\mathrm{a.s.}}\longrightarrow Q$ provided $E (K_{n+1}^2\mid\mathcal
{F}_n )\stackrel{\mathrm{a.s.}}\longrightarrow Q$.

\textit{Almost sure convergence of $E (K_{n+1}^2\mid\mathcal
{F}_n
)$ and $Q>Z^2$ a.s.} Let $G_n=\sum_{x\in S_n}J_n(x)^2$.
Since $M_{n+1}\mid\mathcal{F}_n\sim \operatorname{Be}P(\nu_n)$, then
\[
E \bigl(K_{n+1}^2\mid\mathcal{F}_n \bigr)=
\nu_n(\mathcal{X})+\nu _n(\mathcal{X})^2-\sum
_{x\in\mathcal{X}}\nu_n\{x\}^2=Z_n+Z_n^2-G_n\qquad
\mbox{a.s.;}
\]
see Section~\ref{sh6f}. Thus, since $Z_n\stackrel
{\mathrm{a.s.}}\longrightarrow Z$ and $(G_n)$ is uniformly integrable, it
suffices to prove that $(G_n)$ is a sub-martingale with respect to
$(\mathcal{G}_n)$.

Let us define the random variables $\{T_{n,r}\dvtx n, r\geq1\}$, with
values in $\mathcal{X}\cup\{\infty\}$, as follows. For $n=1$, let
$T_{1,r}=\infty$ for $r>L_1$. If $L_1>0$, define $T_{1,1},\ldots,\break T_{1,L_1}$ to be the dishes tried by customer 1. By induction, at step
$n\geq2$, let
\[
T_{n,r}=T_{n-1,r}\qquad \mbox{for } 1\leq r\leq L_{n-1}\quad
\mbox {and}\quad T_{n,r}=\infty\qquad\mbox{for } r>L_n.
\]
If $L_n>L_{n-1}$, define $T_{n,L_{n-1}+1},\ldots,T_{n,L_n}$ to be the
dishes tried for the first time by customer $n$. Then, $\sigma
(T_{n,r})\subset\mathcal{G}_n$ for all $r\geq1$. Letting $J_n(\infty
)=0$, one also obtains
\[
G_n=\sum_rJ_n
(T_{n,r} )^2.
\]
For fixed $r$, since $\sigma(T_{n,r})\subset\mathcal{G}_n$, it
follows that
\begin{eqnarray*}
&&E \bigl\{J_{n+1}(T_{n,r})\mid\mathcal{G}_n \bigr
\}\\
&&\qquad=E \bigl\{I_{\{
r\leq
L_n\}} J_{n+1}(T_{n,r})\mid
\mathcal{G}_n \bigr\}
\\
&&\qquad=I_{\{r\leq L_n\}} \frac{-\beta+\sum_{i=1}^nR_iM_i\{T_{n,r}\}
+R_{n+1}E \{M_{n+1}\{T_{n,r}\}\mid\mathcal{G}_n \}}{c+\sum_{i=1}^{n+1}R_i}
\\
&&\qquad=I_{\{r\leq L_n\}} \frac{J_n(T_{n,r}) \{c+\sum_{i=1}^nR_i\}
+R_{n+1}J_n(T_{n,r})}{c+\sum_{i=1}^{n+1}R_i}=J_n(T_{n,r})\qquad
\mbox{a.s.}
\end{eqnarray*}
Then
\begin{eqnarray*}
E \{G_{n+1}\mid\mathcal{G}_n \}&=&E \biggl\{\sum
_rJ_{n+1}(T_{n+1,r})^2\Bigm|
\mathcal{G}_n \biggr\}\geq E \biggl\{\sum
_rJ_{n+1}(T_{n,r})^2\Bigm|
\mathcal{G}_n \biggr\}
\\
&=&\sum_rE \bigl\{J_{n+1}(T_{n,r})^2
\mid\mathcal{G}_n \bigr\}\geq \sum_rE
\bigl\{J_{n+1}(T_{n,r})\mid\mathcal{G}_n \bigr
\}^2\\
&=&\sum_rJ_n(T_{n,r})^2=G_n\qquad
\mbox{a.s.}
\end{eqnarray*}
Therefore, $(G_n)$ is a $(\mathcal{G}_n)$-sub-martingale, as required.

From now on, $Q$ denotes a real random variable satisfying
\[
E \bigl(K_{n+1}^2\mid\mathcal{F}_n \bigr)
\stackrel {\mathrm{a.s.}}\longrightarrow Q\quad \mbox{and} \quad\frac{1}{n}\sum
_{i=1}^nK_i^2\stackrel
{\mathrm{a.s.}}\longrightarrow Q.
\]

Let us prove $Q>Z^2$ a.s. Let $Y_n=J_n(T_{n,1})$. Since $(Y_n)$ is a
$[0,1]$-valued sub-martingale with respect to $(\mathcal{G}_n)$, one
obtains $Y_n\stackrel{\mathrm{a.s.}}\longrightarrow Y$ for some random variable~$Y$. Thus
\begin{eqnarray*}
Q-Z^2&=&\lim_n \bigl\{E \bigl(K_{n+1}^2
\mid\mathcal{F}_n \bigr)-Z_n^2 \bigr\}=\lim
_n \{Z_n-G_n \}
\\
&=&\lim_n \biggl\{\Lambda_n+\sum
_rJ_n(T_{n,r}) \bigl(1-J_n(T_{n,r})
\bigr) \biggr\}
\\
&\geq&\lim_n Y_n (1-Y_n)=Y (1-Y)\qquad
\mbox{a.s.}
\end{eqnarray*}
Since $L_n\stackrel{\mathrm{a.s.}}\longrightarrow\infty$ (because of Theorem~\ref
{thmln}), then $T_{n,1}\neq\infty$ eventually a.s. Hence, arguing as in
\cite{AMS} and \cite{MF} (see also Section~4.3 of \cite{BCPR11}) it
can be shown that $Y$ has a diffuse distribution. Therefore, $0<Y<1$
and $Q-Z^2\geq Y (1-Y)>0$ a.s.

We next turn to conditions (jj)--(jjj)--(jv).

\textit{Condition} (jj). Since $E(Z_{n-1}^4)=E \{E(K_n\mid
\mathcal
{F}_{n-1})^4 \}\leq E(K_n^4)$, then
\[
\sup_n E \bigl\{K_n^4+Z_{n-1}^4+
\Lambda_{n-1}^4+N_n^4 \bigr\}\leq2
\sup_n E \bigl\{K_n^4+
\Lambda_{n-1}^4+N_n^4 \bigr\}<\infty.
\]
Therefore,
\begin{eqnarray*}
&&E \Bigl\{ \Bigl(\sup_{n\geq1}\sqrt{n} |Z_n-Z_{n-1}|
\Bigr)^4 \Bigr\} \\
&&\qquad\leq\sum_{n=1}^\infty
n^2E \bigl\{(Z_n-Z_{n-1})^4 \bigr\}
\\
&&\qquad\leq D_1 \sum_{n=1}^\infty
n^2 \biggl\{E \bigl\{ (\Lambda _n-\Lambda _{n-1}
)^4 \bigr\}+\frac{E \{K_n^4+Z_{n-1}^4+\Lambda
_{n-1}^4+N_n^4 \}}{n^4} \biggr\}
\\
&&\qquad\leq D_2 \sum_{n=1}^\infty
\biggl\{\frac{1}{n^{6-4\beta}}+\frac
{1}{n^2} \biggr\}<\infty,
\end{eqnarray*}
where $D_1$ and $D_2$ are suitable constants.

In order to prove (jjj)--(jv), we let
\[
U=\frac{q-r^2}{r^2} \bigl(Q-Z^2 \bigr) \quad\mbox{and}\quad V=
\frac{q}{r^2} \bigl(Q-Z^2 \bigr).
\]

\textit{Condition} (jjj). Let $X_n= \{
K_n-Z_{n-1}+n(Z_{n-1}-Z_n) \}
^2$. On noting that $\sum_nn^2E \{(Z_n-Z_{n-1})^4 \}<\infty$,
as shown in (jj), one obtains $\sum_nn^{-2}E(X_n^2)<\infty$. Thus, by
Lemma~\ref{iam7t}, it suffices to prove $E (X_n\mid\mathcal
{F}_{n-1} )\stackrel{\mathrm{a.s.}}\longrightarrow U$. To this end, we first
note that
\[
E \bigl\{(K_n-Z_{n-1})^2\mid
\mathcal{F}_{n-1} \bigr\}=E \bigl\{ K_n^2\mid
\mathcal{F}_{n-1} \bigr\}-Z_{n-1}^2\stackrel{\mathrm{a.s.}}
\longrightarrow Q-Z^2.
\]
We next prove
\begin{eqnarray*}
&&\mbox{(*)}\quad n^2E \biggl\{\frac{R_n^2(K_n-Z_{n-1})^2}
{ (\sum_{i=1}^nR_i+c )^2}\Bigm|\mathcal{F}_{n-1}
\biggr\}\stackrel {\mathrm{a.s.}}\longrightarrow V; \\
&&\mbox{(**)}\quad n E \biggl\{
\frac
{R_n(K_n-Z_{n-1})^2}{\sum_{i=1}^nR_i+c}\Bigm|\mathcal{F}_{n-1} \biggr\} \stackrel{\mathrm{a.s.}}
\longrightarrow Q-Z^2.
\end{eqnarray*}
In fact,
\begin{eqnarray*}
&&n^2E \biggl\{\frac{R_n^2(K_n-Z_{n-1})^2}{(\sum_{i=1}^nR_i+c)^2}\Bigm| \mathcal{F}_{n-1}
\biggr\}\\
&&\qquad\leq n^2\frac{E \{
R_n^2(K_n-Z_{n-1})^2\mid
\mathcal{F}_{n-1} \}}{ (\sum_{i=1}^{n-1}R_i )^2}
\\
&&\qquad= \biggl(\frac{n}{n-1} \biggr)^2 \frac{E(R_n^2) E \{
(K_n-Z_{n-1})^2\mid\mathcal{F}_{n-1} \}}{ (\overline
{R}_{n-1} )^2}
\stackrel{\mathrm{a.s.}}\longrightarrow\frac{q (Q-Z^2)}{r^2}=V.
\end{eqnarray*}
Since $R_n\leq b$, one also obtains
\[
n^2E \biggl\{\frac{R_n^2(K_n-Z_{n-1})^2}{(\sum_{i=1}^nR_i+c)^2}\Bigm| \mathcal{F}_{n-1}
\biggr\}\geq n^2\frac{E \{
R_n^2(K_n-Z_{n-1})^2\mid
\mathcal{F}_{n-1} \}}{ (\sum_{i=1}^{n-1}R_i+b+c
)^2}\stackrel {\mathrm{a.s.}}\longrightarrow V.
\]
This proves condition (*). Similarly, (**) follows from
\begin{eqnarray*}
n E \biggl\{\frac{R_n(K_n-Z_{n-1})^2}{\sum_{i=1}^nR_i+c}\Bigm|\mathcal {F}_{n-1} \biggr\}& \leq&
\frac{n}{n-1} \frac{E(R_n) E \{(K_n-Z_{n-1})^2\mid\mathcal
{F}_{n-1} \}}{\overline{R}_{n-1}} \\
&\stackrel{\mathrm{a.s.}}\longrightarrow&
Q-Z^2\quad\mbox{and}
\\
 n E \biggl\{\frac{R_n(K_n-Z_{n-1})^2}{\sum_{i=1}^nR_i+c}\Bigm|\mathcal{F}_{n-1} \biggr
\} &\geq& n \frac{E(R_n) E \{(K_n-Z_{n-1})^2\mid\mathcal
{F}_{n-1}
\}}{\sum_{i=1}^{n-1}R_i+b+c}\\
&\stackrel{\mathrm{a.s.}}\longrightarrow& Q-Z^2.
\end{eqnarray*}
Finally, by Lemma~\ref{basic123l} and after some calculations, one obtains
\begin{eqnarray*}
n^2E \bigl\{(Z_{n-1}-Z_n)^2\mid
\mathcal{F}_{n-1} \bigr\}-n^2E \biggl\{ \frac
{R_n^2(K_n-Z_{n-1})^2}{(\sum_{i=1}^nR_i+c)^2}
\Bigm|\mathcal {F}_{n-1} \biggr\}&\stackrel{\mathrm{a.s.}}\longrightarrow&0,
\\
n E \bigl\{(K_n-Z_{n-1}) (Z_{n-1}-Z_n)
\mid\mathcal{F}_{n-1} \bigr\} +n E \biggl\{\frac{R_n(K_n-Z_{n-1})^2}{\sum_{i=1}^nR_i+c}\Bigm|
\mathcal {F}_{n-1} \biggr\}&\stackrel{\mathrm{a.s.}}\longrightarrow&0.
\end{eqnarray*}
Therefore, $n^2E \{(Z_{n-1}-Z_n)^2\mid\mathcal{F}_{n-1} \}
\stackrel{\mathrm{a.s.}}\longrightarrow V$ and
\[
2 n E \bigl\{(K_n-Z_{n-1}) (Z_{n-1}-Z_n)
\mid\mathcal{F}_{n-1} \bigr\} \stackrel{\mathrm{a.s.}}\longrightarrow-2
\bigl(Q-Z^2 \bigr),
\]
which in turn implies $E (X_n\mid\mathcal{F}_{n-1}
)\stackrel
{\mathrm{a.s.}}\longrightarrow V-(Q-Z^2)=U$.

\textit{Condition} (jv). Let $X_n=n^2(Z_n-Z_{n-1})^2$. Since $\sum_nn^2E \{(Z_n-Z_{n-1})^4 \}<\infty$ and $n^2E \{
(Z_{n-1}-Z_n)^2\mid\mathcal{F}_{n-1} \}\stackrel
{\mathrm{a.s.}}\longrightarrow V$, as shown in (jj) and (jjj), Lemma~\ref
{iam7t} yields
\[
n\sum_{k\geq n}(Z_{k-1}-Z_k)^2=n
\sum_{k\geq n}\frac
{X_k}{k^2}\stackrel {\mathrm{a.s.}}
\longrightarrow V.
\]

In view of Theorem~\ref{j2a87p}, this concludes the proof of the first part.

Finally, suppose $R_n=1$ for all $n$. Then, by Lemma~\ref{c4p}, $(M_n)$
is c.i.d. with respect to the filtration
$(\mathcal{G}_n)$. Thus, $(M_n)$ is c.i.d. with respect to
$(\mathcal{F}_n)$ as well, and condition (j) (with
$\mathcal{U}_n=\mathcal{F}_n$) is automatically true. To complete the
proof, it suffices to note that $\beta<1/2$ is only needed in condition
(j). All other points of this proof are valid for each $\beta<1$.
\end{pf*}

\section*{Acknowledgments}
Irene Crimaldi is a member of the Gruppo Nazionale per
l'Analisi Matematica, la Probabilit\`{a} e le loro Applicazioni (GNAMPA)
of the Istituto Nazionale di Alta Matematica (INdAM).
We also thank Paolo
Boldi for suggesting us to investigate IBP-type models, and an
anonymous referee for various helpful remarks.
% imsref loaded by akundreckaite, 2014-03-27 08:19:09
%

%


%suskaldyti doi

\printaddresses

\end{document}